\documentclass[12pt]{article}
\begin{document}
\newtheorem{prop}{}[section]
\newtheorem{defi}[prop]{}
\newtheorem{lemma}[prop]{}
\newtheorem{rema}[prop]{}
\def\s{\scriptstyle}
\def\ss{\scriptscriptstyle}
\def\dd{\displaystyle}
\def\Hn{ H^{n}(\reali^d, \complessi) }
\def\Ld{L^{2}(\reali^d, \complessi)}
\def\Lpi{L^{p}(\reali^d, \complessi)}
\def\Lq{L^{q}(\reali^d, \complessi)}
\def\Lr{L^{r}(\reali^d, \complessi)}
\def\k{\mbox{{\tt k}}}
\def\x{\mbox{{\tt x}}}
\def\lgraffa{ \mbox{\Large $\{$ } \hskip -0.2cm}
\def\rgraffa{ \mbox{\Large $\}$ } \hskip -0.2cm}
\def\sc{ {\scriptstyle{\bullet} }}
\def\la{\lambda}
\def\complessi{{\textbf C}}
\def\reali{{\textbf R}}
\def\interi{{\textbf Z}}
\def\naturali{{\textbf N}}
\def\FF{{\mathcal F}}
\def\SS{{\mathcal S}}
\def\circa{\thickapprox}
\def\vain{\rightarrow}
\def\parn{\par \noindent}
\def\salto{\vskip 0.2truecm \noindent}
\def\spazio{\vskip 0.5truecm \noindent} 
\def\vs1{\vskip 1cm \noindent}
\def\fine{\hfill $\diamond$}
\newcommand{\rref}[1]{(\ref{#1})}
\def\beq{\begin{equation}}
\def\feq{\end{equation}}
\def\beqq{\begin{eqnarray}}
\def\feqq{\end{eqnarray}}
\def\barray{\begin{array}}
\def\farray{\end{array}}
\makeatletter
\@addtoreset{equation}{section}
\renewcommand{\theequation}{\thesection.\arabic{equation}}
\makeatother
\begin{center}
{\huge On the constants for some Sobolev imbeddings} 
\end{center}
\vspace{0.4truecm}
\begin{center}
{\large
Carlo Morosi${}^1$, Livio Pizzocchero${}^2$} \\
\vspace{0.5truecm}
{\footnotesize
${}^1$ Dipartimento di Matematica, Politecnico di
Milano, \\ P.za L. da Vinci 32, I-20133 Milano, Italy \\
e--mail: carmor@mate.polimi.it \\
${}^2$ Dipartimento di Matematica, Universit\`a di Milano\\
Via C. Saldini 50, I-20133 Milano, Italy\\
and Istituto Nazionale di Fisica Nucleare, Sezione di Milano, Italy \\
e--mail: livio.pizzocchero@mat.unimi.it \\ 
Work partly supported by MURST and Indam, Gruppo Nazionale 
per la Fisica Matematica. }
\end{center}
\begin{abstract} \noindent
We consider the imbedding inequality $\|~ \|_{L^r(\reali^d)}
\leq S_{r,n, d}~ \|~\|_{H^{n}(\reali^d)}$; $H^{n}(\reali^d)$ is the 
Sobolev space (or Bessel potential space) of $L^2$ type
and (integer or fractional) order $n$. We write down upper bounds 
for the constants $S_{r, n, d}$, using an argument 
previously applied in the literature in particular cases. 
We prove that the upper bounds computed 
in this way are in fact the sharp constants if
($r=2$ or) $n > d/2$, $r=\infty$, and exhibit the maximising functions. 
Furthermore, using convenient trial functions, we derive 
lower bounds on $S_{r,n,d}$ for $n > d/2$, $2 < r < \infty$; in many
cases
these are close to the previous upper bounds, as illustrated by a 
number of examples, thus characterizing the sharp constants  
with little uncertainty. 
\end{abstract}
\vskip 0.2cm \noindent
\textbf{To appear in the "Journal of inequalities and applications".}
\section{Introduction and preliminaries.} 
\label{intro}
The imbedding inequality of $\Hn$ into $L^r(\reali^d, \complessi)$ is a 
classical topic, and several approaches has been developed to derive 
upper bounds on the sharp imbedding constants $S_{r,n,d}$. A simple
method, 
based on the Hausdorff-Young and H\"older inequalities, has been 
employed in the literature for special choices of $r,n,d$, as 
indicated in the references at the end of Sect.\ref{up}. 
Little seems to have been done 
to test reliability of the upper bounds derived in this 
way (i.e., their precision in approximating the unknown 
sharp constants). \parn
This paper is a contribution to the understanding of the
Hausdorff-Young-H\"older (HYH) upper bounds, and aims to show their
reliability for $n > d/2$. This case is interesting for a number 
of reasons, including application to PDE's; its main feature is that
the $H^n$ norm controls the $L^r$
norms of \textsl{all} orders $r \geq 2$, up to $r=\infty$. \parn
The paper is organized as follows. First of all, 
in Sect.\ref{up} we write the general expression of the HYH  
upper bounds $S_{r,n,d} \leq S^{+}_{r,n,d}$ (containing all special
cases 
of our knowledge in the literature). 
In Sect.\ref{low} we show that the upper 
bounds $S^{+}_{r,n,d}$ are in fact the sharp constants if 
($r=2$, $n$ arbitrary or) $n > d/2$, $r=\infty$, 
and exhibit the maximising functions; next, we assume $n > d/2$
and inserting a one parameter family
of trial functions in the imbedding inequality, we derive lower bounds
$S_{r,n,d} \geq S^{-}_{r,n,d}$ for arbitrary $r \in (2, \infty)$. 
In Sect.\ref{ese} we report numerical values of $S^{\pm}_{r,n,d}$ 
for representative choices of $n, d$ and a wide range of $r$ values; 
in all the examples the relative uncertainty on the sharp imbedding 
constants, i.e. the ratio
$(S^{+}_{r,n,d} - S^{-}_{r,n,d})/S^{-}_{r,n,d}$, is found to be $\ll 1$. 
\vskip 0.1cm \noindent
\textbf{Notations for Fourier transform and $H^{n}$ spaces.} 
Throughout this paper, $d \in \naturali \setminus \{ 0 \}$ is a 
fixed space dimension; the running variable in $\reali^d$ is 
$x = (x_1, ..., x_d)$, and $k = (k_1, ..., k_d)$ when using the 
Fourier transform. We write $| \x |$ 
for the function $(x_1, ..., x_d) \mapsto \sqrt{{x_1}^2 + ... +
{x_d}^2}$, and
intend $| \k |$ similarly.
We denote with 
$\FF, \FF^{-1} : \SS'(\reali^d, \complessi) \vain \SS'(\reali^d, 
\complessi)$
the Fourier transform of tempered distributions 
and its inverse, choosing normalizations so that
(for $f$ in $L^1(\reali^d, \complessi)$~) it is
$\FF f(k)= (2 \pi)^{-d/2} \times$ 
$\times \int_{\reali^d} d x~e^{-i k \sc x}
f(x)$.
The restriction of $\FF$ to $\Ld$, with the standard
inner product and the associated norm $\|~\|_{L^2}$, is 
a Hilbertian isomorphism. \parn
For real $n \geq 0$, let us introduce the operators 
\beq \SS'(\reali^d, \complessi) \vain \SS'(\reali^d, \complessi)~, 
\qquad g \mapsto \sqrt{1 - \Delta}^{~\pm n}~ g 
:= \FF^{-1} \left( \sqrt{1 + | \k |^2}^{~\pm n} \FF g \right)
\label{lap} \feq
(in case of integer, even exponent $n$, we have a power of 
$1$ minus the distributional Laplacian $\Delta$, in the elementary
sense). The $n$-th order Sobolev
(or Bessel potential \cite{Smi}) space of $L^2$ type and its norm are 
$$ \Hn := \lgraffa f \in \SS'(\reali^d, \complessi)~\Big\vert~
\sqrt{1 - \Delta}^{~n} f \in \Ld~ \rgraffa= $$
$$ = \lgraffa 
\sqrt{1 - \Delta}^{~-n}~ u~\Big\vert~u \in \Ld~\rgraffa= $$
\beq = \lgraffa f \in \SS'(\reali^d, \complessi)~\Big \vert~
\sqrt{1 + | \k |^2}^{~n}  \FF f \in \Ld \rgraffa~, \label{incid} \feq
\beq \| f \|_{H^n} := \| \sqrt{1 - \Delta}^{~n}~ f \|_{L^2} = 
\|~ \sqrt{1 + | \k |^2}^{~n}~\FF f~ \|_{L^2}~. \label{repfur} \feq
Of course, if $n \leq n'$, it is 
$H^{n'}(\reali^d, \complessi) \subset \Hn$ and $\|~\|_{H^n} \leq
\|~\|_{H^{n'}}$; 
also, $H^{0} = L^2$. 
\vskip 0.1cm \noindent
\textbf{Connection with Bessel functions.} For $\nu > 0$, and 
in the limit case zero, let us put, respectively,  
\beq G_{\nu, d} :=  
\FF^{-1}\left({1 \over \sqrt{1 + | \k |^2}^{~\nu}}\right) = 
{| \x |^{\nu/2 - d/2} \over 2^{\nu/2 - 1} \Gamma(\nu/2)}~
K_{\nu/2 - d/2}(| \x |); \label{gnud} \feq
$$ G_{0, d} := \FF^{-1} (1)= (2 \pi)^{d/2}~\delta. $$
Here, $\Gamma$ is the factorial function; $K_{(~)}$ are the
modified Bessel function of the third kind, or Macdonald functions, see 
e.g. \cite{Wat}; $\delta$ is the Dirac distribution. 
The expression of $G_{\nu, d}$ via a Macdonald function  
\cite{Smi} comes from the known computational rule for
the Fourier transforms of radially symmetric functions \cite{Boc}. 
With the above ingredients, we obtain another representation of $H^n$
spaces 
\cite{Smi}; in fact, explicitating $\sqrt{1 - \Delta}^{~-n}~ u$ in 
Eq.\rref{incid} and recalling
that $\FF^{-1}$ sends pointwise product into $(2 \pi)^{-d/2}$ times
the convolution product $\ast$, we see that 
\beq \Hn = \lgraffa {1 \over (2 \pi)^{d/2}}~G_{n, d} \ast 
u~\Big\vert~u \in \Ld~\rgraffa~ \feq
for each $n \geq 0$. All this is standard; in this paper we will show
that, for 
$n > d/2$, the function 
$G_{2 n, d}$ also plays a relevant role for $\Hn$, being an 
element of this space and appearing to be a
maximiser for the inequality $\| ~\|_{L^{\infty}} \leq 
\mbox{const}~ \| ~\|_{H^n}$. Incidentally we note that 
(for all $n \geq 0$) the relation $(1 + | \k |^2)^{-n} = 
\sqrt{1 + | \k |^2}^{~-n}~ \sqrt{1 + | \k |^2}^{~-n}$ gives, after 
application of $\FF^{-1}$, 
$ G_{2 n, d} = (2 \pi)^{-d/2}~G_{n, d} \ast G_{n, d}$. \parn
For future conveniency, let us recall a case in which the 
expression of $G_{\nu, d}$ simply involves 
an exponential $\times$ a polynomial in $| \x |$. This occurs if $\nu/2
- 
d/2 = m + 1/2$, with $m$ a nonnegative integer: 
in fact, it is well known \cite{Wat} that
\beq \rho^{m+1/2}~K_{m+1/2}(\rho) = 
\sqrt{{\pi \over 2}}~e^{-\rho}~ 
\sum_{i=0}^m {( 2 m - i)! \over i! (m - i)!}~ {\rho^{i} \over 2^{m-i}}
\qquad (m \in \naturali,~\rho \in \reali)~. 
\label{cita} \feq
\vskip 0.4cm \noindent
\section{HYH upper bounds for the imbedding constants.} 
\label{up}
It is known \cite{Smi} \cite{Ada} that $\Hn$ is continuously imbedded
into
$L^{r}(\reali^d, \complessi)$ if 
$0 \leq n < d/2$, $2 \leq r \leq d/(d/2 - n)$ or
$n = d/2$, $2 \leq r < \infty$ or $n > d/2$, $2 \leq r \leq \infty$.
We are interested
in the sharp imbedding constants
\beq S_{r,n, d} := \mbox{Inf}~\{ S \geq 0~|~\mbox{ 
$\| f \|_{L^r} \leq S~\| f \|_{H^n}$~for 
all $f \in \Hn$}~\}~. \label{ineq} \feq
Let us derive general upper bounds on the above constants, with the
HYH method mentioned in the Introduction; this result will be expressed
in 
terms of the functions $\Gamma$ and $E$, the latter being defined by 
\beq E(s) := s^s \qquad \mbox{for $s \in (0, + \infty)$}~, 
\qquad E(0) := \mbox{lim}_{s \vain 0^{+}} E(s) = 1~. \feq
\begin{prop}
\label{sob}
\textbf{Proposition.} Let $n=0$, $r=2$ or
$0 < n < d/2$, $2 \leq r < d/(d/2 - n)$ or
$n = d/2$, $2 \leq r < \infty$ or $n > d/2$, $2 \leq r \leq \infty$.
Then $S_{r,n,d} \leq S^{+}_{r,n,d}$, where
$$ S^{+}_{r, n, d} :=
{1 \over (4 \pi)^{d /4 - d/(2 r)} }~
\left( { \Gamma \left( {\displaystyle{ {n \over 1 - 2/r} - {d \over 2}
}} \right) 
\over \Gamma
\left( {\displaystyle{ {n \over 1 - 2/r} }} \right) } 
\right)^{1/2 - 1/r} \times $$
\beq \times \left({E(1/r) \over E(1 - 1/r)} \right)^{d/2}
\mbox{if~ $r \neq 2, \infty$}~, \label{srnd} \feq
\beq S^{+}_{2, n, d} :=  1~, 
\qquad 
S^{+}_{\infty,n,d} :=  
{1 \over (4 \pi)^{d /4} }~
\left( { \Gamma( n  - d/2 ) \over \Gamma(n) } 
\right)^{1/2}~. 
 \label{srid} \feq
\end{prop}
\textbf{Proof.} Of course, it amounts to showing that 
$\| f \|_{L^r} \leq  S^{+}_{r,n,d}~ \| f \|_{H^n}$
for all $f \in \Hn$.
For $r=2$ and any $n$ this follows trivially, 
because $\| f \|_{L^2} = \| f \|_{H^0}$ $\leq 1 \times \| f \|_{H^n}$.
\parn
>From now on we assume $r \neq 2$ (intending $1/r := 0$
if $r = \infty$); $p$, $s$ are such that
\beq {1 \over r} + {1 \over p} = 1~; \qquad 
{1 \over s} + {1 \over 2} = {1 \over p}~, \quad 
\mbox{i.e.}~, \qquad s= {2 \over 1 - 2/r}~. \label{defes} 
\feq
Let $f \in \Hn$. Then, the Hausdorff-Young
inequality for $\FF$ and  the (generalized) H\"older's inequality for
$ \FF f = \sqrt{1 + | \k |^2}^{~-n} 
\left( {\sqrt{1 + | \k |^2}}^{~n}~\FF f \right)$
give
\beq \| f \|_{L^r} \leq C_{r, d}~ \| \FF f \|_{L^p}~, \qquad 
C_{r, d} := {1 \over (2 \pi)^{d/2 - d/r}}~
\left( {E(1/r) \over E(1 - 1/r)} \right)^{d/2}~, \label{crd} \feq
$$ \| \FF f \|_{L^p} \leq \| {1 \over \sqrt{1 + | \k |^2}^{~n}}
\|_{L^s}~
\| \sqrt{1 + | \k |^2}^{~n} \FF f \|_{L^2} = $$ 
\beq = \left( \int_{\reali^d} 
{d k \over \sqrt{1 + | k |^2}^{~n s} }  \right)^{1/s}
\| f \|_{H^n} \label{noff} \feq
($C_{r,d}$ is the sharp Hausdorff-Young constant: 
see \cite{Lie}, \cite{Lie2} Chapter 5 and references therein. 
Our expression for $C_{r, d}$ differs by a factor from the one
in \cite{Lie2} due to another
normalization for the Fourier transform). \parn
Of course, statements \rref{crd} \rref{noff} are meaningful if
the integral in Eq.\rref{noff} converges; in fact this is the case, 
because the definition of $s$ and the assumptions on $r,n,d$ imply  
$n s > d$.
Summing up, we have
\beq \| f \|_{L^r} \leq {1 \over (2 \pi)^{d/2 - d/r}}~
\left( {E(1/r) \over E(1 - 1/r)} \right)^{d/2}~
\left( \int_{\reali^d} 
{d k \over \sqrt{1 + | k |^2}^{~n s} }  \right)^{1/s}~
\| f \|_{H^n}~,  
\label{give} \feq
with $s$ as in \rref{defes}.
Now, the thesis is proved if we show that
\beq \mbox{constant in Eq.\rref{give}} = S^{+}_{r,n,d}~; \label{word}
\feq
to check this, it suffices to write
\beq \int_{\reali^d} 
{d k \over \sqrt{1 + | k |^2}^{~n s}} = 
{ 2 \pi^{d/2} \over \Gamma(d/2)}~\int_{0}^{+\infty} d \xi 
{ \xi^{d-1}\over \sqrt{1 + \xi^2}^{~n s} } = 
\pi^{d/2}~{\Gamma(n s/2 - d/2) \over \Gamma(n s/2) }~  \label{integ}
\feq
and to explicitate $s$. \fine 
\vskip 0.1cm \noindent
\textbf{Remarks.} i) Let us indicate the special cases of our knowledge,
in 
which some HYH upper bounds $S^{+}_{r,n,d}$ 
have been previously given in the literature. 
Reference \cite{Lie2} derives these bounds  
for $d=1$, $n=1/2$, $d=2$, $n=1$ and $2 \leq r < \infty$ 
(with a misprint). The inequality in \cite{Ree}, page 55 is 
strictly related to the case $n=2$, $d \geq 4$, $2 \leq r < d/(d/2 -
2)$. 
The upper bound $S^{+}_{\infty,n,d}$ is given for
arbitrary $n > d/2$ by many authors, see e.g. \cite{Miz} \cite{Zei}.
\parn
To our knowledge, little was done to discuss reliability of the HYH 
upper bounds; the next two sections will be devoted to this topic, in 
the $n > d/2$ case. First of all, we will emphasize that 
$S^{+}_{\infty,n,d}$ is in fact the \textsl{sharp} imbedding 
constant for any $n > d/2$ 
(this is shown in \cite{Lie2} for $d=1$, $n=1$ only, 
with an \textsl{ad hoc} technique). $S^{+}_{2,n,d}$ is also the sharp
constant
(for any $n$), by an obvious argument; our analysis 
will show that, for $n > d/2$ and intermediate values $2 < r < \infty$, 
$S^{+}_{r,n,d}$ gives a generally good approximation of the sharp
constant. 
\parn
ii) Discussing reliability of the bounds $S^{+}_{r,n,d}$ for 
$n \leq d/2$ would require a separate 
analysis, which is outside the purposes of this paper; let us only
present a few
comments. \parn
The upper bound $S^{+}_{r,n,d}$ is certainly far from the sharp constant 
for $0 < n < d/2$ and $r$ close to $d/(d/2 -n)$: 
note that $S^{+}_{r,n,d}$ diverges for 
$r$ approaching this limit, 
in spite of the validity of the imbedding inequality even at
the limit value. 
As a matter of fact other approaches, not using 
the HYH scheme, are more suitable to estimate the imbedding constants
if  $0 < n < d/2$, $r \simeq d/(d/2 -n)$. 
We refer, in particular, to methods based on the
Hardy-Littlewood-Sobolev
inequality \cite{Miz}: the sharp constants for that inequality were
found 
variationally in \cite{Lih}. Let us also mention 
the papers \cite{Aub}, prior to 
\cite{Lie}, and \cite{Wan}; the inequalities considered therein,
for which the sharp constants were determined, are 
strictly related to the limit case $r=d/(d/2 -n)$ with $n=1$ and $2$, 
respectively. \parn
The HYH upper bounds $S^{+}_{r,n,d}$ might be close to the sharp 
imbedding constants $S_{r,n,d}$
in the critical case $n=d/2$, but this topic will not be 
discussed in the sequel. 
\vskip 0.4cm \noindent
\section{Cases where $\mbox{\boldmath $S^{+}_{r,n,d}$}$
is the sharp constant. Lower
bounds on the sharp constants for $\mbox{\boldmath $n>d/2$}$ and 
arbitrary $\mbox{\boldmath $r$}$. } 
\label{low} 
Let us begin with the aforementioned statement that
\begin{prop}
\label{fore} 
\textbf{Proposition.} $S^{+}_{r,n,d}$ is the sharp imbedding constant
if $n \geq 0$, $r=2$ or $n > d/2$, $r=\infty$. In fact: \parn
\hskip 0.3cm
i) for any  $n \geq 0$ and nonzero $f \in \Hn$, it is
\beq \lim_{\lambda \mapsto 0^{+}} 
{\| f^{(\lambda)} \|_{L^2} \over \| f^{(\lambda)} \|_{H^n} } = 1 = 
S^{+}_{2,n,d}~, \quad \mbox{where $f^{(\la)}(x) :=  f(\la x)$ 
for $x$ in $\reali^d$, $\lambda \in (0, +\infty)$}~.\label{ovvia} \feq
\beq ii) \quad 
\| f \|_{L^{\infty}} = S^{+}_{\infty, n, d}~ \| f \|_{H^n} 
\quad \mbox{for $n > d/2$ and $\displaystyle{ 
f := \FF^{-1} \left( {1 \over (1 + | \k |^2)^{n}~} \right) = 
G_{2 n, d}}$}~.  \qquad \label{choice} \feq
\end{prop} 
\textbf{Proof.} 
i) Given any $f \in \Hn$, define $f^{(\la)}$ as above; by elementary 
rescaling of the integration variables, we find
\beq \Big(\FF f^{(\la)}\Big)(k) = 
{1 \over \la^d} (\FF f)({k \over \la}) \qquad 
\mbox{for $k \in \reali^d$}~;
\label{resc} \feq
\beq \| f^{(\la)} \|_{H^n} = {1 \over \la^d}~\sqrt{ \int_{\reali^d} d k~
(1 + | k |^2)^n | \FF f({k \over \la}) |^2} = 
{1 \over \la^{d/2}}~\sqrt{ \int_{\reali^d} d h~
(1 + \la^2 | h |^2)^n | \FF f(h) |^2}~; \feq
\beq \| f^{(\la)} \|_{H^n} \barray{ccc} ~ \\ 
\sim \\ \scriptstyle{\la \vain 0^{+}} \farray 
{1 \over \la^{d/2}}~\sqrt{ \int_{\reali^d} d h~
| \FF f(h) |^2} = {1 \over \la^{d/2}} \| f \|_{L^2} = 
\| f^{(\la)} \|_{L^2}~. \label{den} \feq
\parn
ii) Let $n > d/2$; then 
$1/(1 + | \k |^2)^{n} \in L^1(\reali^d, \complessi)$, so 
$f$ in Eq. \rref{choice} is continuous and bounded. 
For all $x \in \reali^d$ (and for the everywhere continuous
representative 
of $f$) it is 
\beq f(x) =
{1 \over (2 \pi)^{d/2} } \int_{\reali^d} d k~ 
{e^{i k \sc x} \over (1 + | k |^2)^{n}~}~, \quad
| f(x) | \leq \int_{\reali^d} d k~ {1\over (1 + | k |^2)^{n}~} = f(0),
\feq
so that
\beq \| f \|_{L^{\infty}} = f(0) = 
{1 \over (2 \pi)^{d/2} } \int_{\reali^d} d k~{1 \over (1 + | k
|^2)^{n}}~. 
\feq
Also, it is $f \in \Hn$ and
\beq \| f \|_{H^n} = \sqrt{ 
\int_{\reali^d} d k~ (1 + | k |^2)^{n}~| \FF f(k) |^2 } = 
\sqrt{ \int_{\reali^d} d k~{1 \over (1 + | k |^2)^{n}}  }~. \feq
The last two equations give
\beq {\| f \|_{L^{\infty}} \over \| f \|_{H^n}} = 
{1 \over (2 \pi)^{d/2} } 
\sqrt{\int_{\reali^d} d k~{1 \over (1 + | k |^2)^{n}}}~,  \feq
and by comparison with Eq.s \rref{give} \rref{word} we see that 
the above ratio is just $S^{+}_{\infty,n,d}$. 
\fine 
\parn
As an example, let us write down the maximising function $f = G_{2 n,
d}$
of item ii) when $n = d/2 + 1/2$ or $n=d/2 + 1$. According to 
Eq.s \rref{gnud} \rref{cita}, we have 
$$ G_{2(d/2 + 1/2), d} = {| \x |^{1/2} K_{1/2}(|\x|) \over 2^{d/2 - 1/2} 
\Gamma(d/2 + 1/2)}  = {\sqrt{\pi}~ e^{-| \x |} 
\over 2^{d/2} \Gamma(d/2 + 1/2)};~ $$
\beq G_{2(d/2 + 1), d} = {| \x |~ K_{1}(|\x|) \over 2^{d/2} 
\Gamma(d/2 + 1)} . \feq
>From now on $n > d/2$; we 
attack the problem of finding lower bounds on $S_{r,n,d}$ for
$2 < r < \infty$. To obtain them, one can insert 
into the imbedding inequality \rref{ineq} a trial function; 
the previous considerations suggest to employ the 
one parameter family of rescaled
functions 
\beq G^{(\la)}_{2 n, d}(x) := G_{2 n, d}(\la x) \quad \qquad 
(~\la \in (0, \infty)~)~.\feq
Of course, the sharp constant satisfies
\beq S_{r, n, d} \geq \mbox{Sup}_{\lambda > 0}~ 
{ \| G^{(\la)}_{2 n, d} \|_{L^r} \over
\| G^{(\la)}_{2 n, d} \|_{H^n} }~;  \label{belo} \feq
one should expect the above supremum to be attained
for $\lambda \simeq 0$ if $r \simeq 2$, and for $\lambda 
\simeq 1$ if $r$ is large. Evaluation of the above 
ratio of norms leads to the following 
\begin{prop}
\label{disot}
\textbf{Proposition.} For $n > d/2$, $2 < r < \infty$ it is 
$S_{r, n, d} \geq S^{-}_{r,n,d}$,
where
\beq S^{-}_{r,n,d} := \left({\Gamma(d/2) \over 2 \pi^{d/2} }\right)^{1/2
- 1/r}~
{ I_{r, n, d}^{1/r} \over 2^{n-1} \Gamma(n) \sqrt{ \Phi_{r,n,d} } }~, 
\label{smrnd} \feq
\beq I_{r,n,d} := \int_{0}^{+\infty} 
d t~ t^{d - 1} \left( t^{n-d/2} K_{n - d/2}(t)\right)^r~, \label{irnd}
\feq
\beq \Phi_{r,n,d} := {\rm{Inf}}_{\lambda > 0}~\varphi_{r,n,d}(\la),~ 
\quad \varphi_{r,n,d}(\la) := {1 \over \la^{d - 2 d/r}}
\int_{0}^{+\infty} d s~
s^{d-1} {(1 + \la^2 s^2)^n \over (1 + s^2)^{2 n}} \label{defirnd}~. \feq
\end{prop}
\textbf{Proof.} From the explicit expression \rref{gnud}, it follows 
(using the variable $t = \la | x |$)
\beq \| G^{(\la)}_{2 n,d} \|^{r}_{L^r} = 
{2 \pi^{d/2} \over \Gamma(d/2)}~
{1 \over 2^{r(n-1)} \Gamma(n)^r~\lambda^d} 
\int_{0}^{+\infty} 
d t~ t^{d - 1} \left( t^{n-d/2} K_{n - d/2}(t) \right)^r~.
\label{l1} \feq
By \rref{resc} with $f=G_{2 n, d}$, it is
$(\FF G^{(\la)}_{2 n, d})(k) = $ $\la^{-d}~(1 + | k |^2/\la^2)^{-n}$,
whence
(using the variable $s = | k |/\la$)
$$ \| G^{(\la)}_{2 n, d} \|_{H^n}^2 = {1 \over \la^{2 d}} 
\int_{\reali^d} d k {(1 + | k |^2)^n \over (1 + | k |^2/\la^2)^{2 n}} =
$$
\beq = {2 \pi^{d/2} \over \Gamma(d/2) \la^{d}} \int_{0}^{+\infty} 
d s~ s^{d-1} {(1 + \la^2 s^2)^n \over (1 + s^2)^{2 n}}~. \label{l2} \feq
Eq.s \rref{l1} \rref{l2} imply 
\beq  { \| G^{(\la)}_{2 n, d} \|_{L^r} \over
\| G^{(\la)}_{2 n, d} \|_{H^n} }
= \left({\Gamma(d/2) \over 2 \pi^{d/2}}\right)^{1/2 - 1/r}
{ I_{r, n, d}^{1/r} \over 2^{n-1} \Gamma(n) \sqrt{ \varphi_{r,n,d}(\la)
} }~, 
\feq
and \rref{belo} yields the thesis. \fine 
\vskip 0.1cm \noindent
\textbf{Remarks.} i) For $n$ integer, the integral in the definition 
\rref{defirnd} of $\varphi_{r,n,d}$ is readily computed  
expanding $(1 + \la^2 s^2)^n$ with the binomial formula,
and integrating term by term. The integral of each term is 
expressible via the Beta 
function $B(z, w) = \Gamma(z) \Gamma(w)/\Gamma(z + w)$, the 
final result being
\beq \varphi_{r,n,d}(\la) = {1 \over 2 \la^{d - 2 d/r}}~\sum_{\ell=0}^n~
\left( \barray{c} n \\ \ell \farray \right)
B(\ell + {d \over 2}, 2 n - {d \over 2} - \ell)~\la^{2 \ell} 
\qquad (n \in \naturali)~. \label{remar} \feq
For arbitrary, possibly noninteger $n$, 
the integral in \rref{defirnd} can be expressed in terms of 
the Gauss hypergeometric function $F= {~}_{2} F_{1}$, and 
the conclusion is
\beq \varphi_{r,n,d}(\la) = {1 \over 2 \la^{d - 2 d/r}}~\Big( 
B(2 n - {d \over 2}, {d \over 2})~ F({d \over 2}, -n, 1 + 
{d \over 2} - 2 n; \la^2) + 
\label{sum} \feq
$$ + \la^{4 n - d}~B(n - {d \over 2}, {d \over 2} - 2 n)~ 
F(2 n, n - {d \over 2}, 1 - {d \over 2} + 2 n; \la^2)~\Big) $$
(in the singular cases $2 n - d/2 - 1 \in \naturali$, 
the first hypergeometric in \rref{sum} 
must be appropriately intended, as a limit from nonsingular values).
\parn
ii) Concerning $I_{r,n,d}$, there is one case in which the integral 
\rref{irnd} is elementary, namely $n = d/2 + 1/2$. In fact, this case 
involves the function $t^{1/2}~ K_{1/2}(t) =$ $\sqrt{\pi/ 2}~e^{-t}$,
so that 
\beq I_{r,d/2 + 1/2,d} = \left({\pi\over 2}\right)^{r/2}~
\int_{0}^{+\infty} d t~t^{d -1} e^{-r t} = 
\left( {\pi \over 2} \right)^{r/2} {\Gamma(d) \over r^d} 
~. \label{elem} \feq
More generally, if $n = d/2 + m + 1/2$, $m \in \naturali$, the
integral defining $I_{r,n,d}$ involves the function 
$t^{m+1/2}~K_{m+1/2}(t)$, which has the elementary expression
\rref{cita};
for $n$ as above and $r$ integer, expanding the power 
$\left(t^{m+1/2}~K_{m+1/2}(t)\right)^r$ 
we can reduce $I_{r,n,d}$ to a linear 
combination of integrals of the type $\int_{0}^{+\infty} t^{\alpha}
e^{-r t} = 
\Gamma(\alpha+1)/r^{\alpha+1}$. 
In other cases, $I_{r,n,d}$ can be evaluated numerically.
\vskip 0.4cm \noindent
\section{Examples.} 
\label{ese}
We present four examples A) B) C) D), each one corresponding to fixed
values of 
$(n, d)$ with $n > d/2$, and 
$r$ ranging freely. Of course, in all these cases the analytical
expression \rref{srnd} of $S^{+}_{r,n,d}$ is available; 
the expressions of the lower bounds 
$S^{-}_{r,n,d}$ are simple in 
examples A) D) and more complicated in examples B) C), where 
the integral $I_{r,n,d}$ is not expressed in terms of elementary
functions,
for arbitrary $r$.
\parn
Each example is concluded by a table of numerical values of 
$S^{\pm}_{r,n,d}$ (computed with the MATHEMATICA package), which 
are seen to be fairly close; the relative uncertainty
$(S^{+}_{r,n,d} - S^{-}_{r,n,d})/S^{-}_{r,n,d}$ is also evaluated.
In cases A) C) D) the space dimension is $d=1,2,3,$ 
respectively, and we take for 
$n$ the smallest integer $>d/2$: this choice of $n$ is the most
interesting
in many applications to PDE's. 
In case B) where $n$ is larger, the uncertainty is even smaller.
Whenever we give numerical values, we 
round from above the digits of $S^{+}_{r,n,d}$, 
and from below the digits of $S^{-}_{r,n,d}$.
\vskip 0.1cm \noindent
\textbf{A) Case $\mbox{\boldmath $n=1, d=1$}$}~. Eq.s \rref{srnd} 
\rref{srid} give $S^{+}_{r,1,1}$ for all $r \in [2, \infty]$; 
the values at the extremes are
\beq S^{+}_{2, 1, 1} = 1~, \qquad S^{+}_{\infty, 1, 1} = 1/\sqrt{2} 
\simeq 0.7072 \feq 
(coinciding with the sharp imbedding constants due to Prop.\ref{fore}). 
Let us pass to the lower bounds. The function $\varphi_{r,1,1}$ is
given by \rref{remar} and attains its minimum at a point 
$\la = \la_{r,1,1}$; the integral $I_{r,1,1}$ is provided by 
\rref{elem}, and these objects must be inserted into \rref{smrnd}.
Explicitly,
\beq \varphi_{r,1,1}(\la) = {\pi (\la^2 + 1)
\over 4~ \la^{1 - 2/r} }~, \quad \la_{r,1,1} = \sqrt{{1-2/r \over 1 +
2/r}}~; 
\qquad I_{r,1,1} = \left( {\pi \over 2} \right)^{r/2}~ {1 \over r} 
~; \feq
\beq S^{-}_{r,1,1} = { I_{r,1,1}^{1/r} \over 2^{1/2 - 1/r} 
\sqrt{\varphi_{r,1,1}(\la_{r,1,1})} } = 
{ E(1/r) \over 2^{1/2 - 1/r}}~
E\left(1 + {2 \over r} \right)^{1/4} E \left(1 - {2 \over r}
\right)^{1/4} 
~.\label{sm11} \feq
Computing numerically the bounds \rref{srnd} \rref{sm11} for many 
values of $r \in (2, +\infty)$, we always found 
$(S^{+}_{r,1,1} - S^{-}_{r,1,1})/S^{-}_{r,1,1} < 0.05$,
the maximum of this relative uncertainty being attained 
for $r \simeq 6$. Here are some 
numerical values:
\beq \begin{tabular}{c|c|c|c|c|c|c} $r$ & 2.2 & 3 & 4 & 6 & 50 & 1000 \\
\hline
$S^{+}_{r,1,1}$ & 0.8832 & 0.7212 & 0.6624 & 0.6345 & 0.6782 & 0.7046 \\
\hline
$S^{-}_{r,1,1}$ & 0.8730 & 0.6973 & 0.6347 & 0.6057 & 0.6632 & 0.7027  
\end{tabular}~~~~. \feq
\vskip 0.1cm \noindent
\textbf{B) Case $\mbox{\boldmath $n=3, d=1$}$}~. Eq.s \rref{srnd} 
\rref{srid} give $S^{+}_{r,3,1}$ for all $r$; in particular
\beq S^{+}_{2, 3, 1} = 1~, \qquad S^{+}_{\infty, 3, 1} = 
\sqrt{3}/4 \simeq 0.4331~.\feq 
We pass to the lower bounds. Eq.s \rref{remar} \rref{irnd} \rref{cita} 
give
$$ \varphi_{r,3,1}(\la) = {3 \pi (\la^6 + 3 \la^4 + 7 \la^2 + 21)
\over 512~ \la^{1 - 2/r}}~; $$
\beq I_{r,3,1} = \left( {\pi \over 2} \right)^{r/2} \int_{0}^{+\infty} d
t~
(t^2 + 3 t + 3)^r e^{-r t}~.
\feq
The minimum point $\la_{r,3,1}$ of $\varphi_{r,3,1}$ is the positive 
solution of the equation
\beq (5 + {2 \over r}) \la^6 + (9 + {6 \over r}) \la^4 + 
(7 + {14 \over r}) \la^2 - (21 - {42 \over r}) = 0~; \feq
the integral $I_{r,3,1}$ 
can be computed analytically for integer $r$, and numerically 
otherwise. The final lower bounds, and some numerical values for 
them and for the upper bounds \rref{srnd} are
\beq S^{-}_{r,3,1} = { I_{r,3,1}^{1/r} \over 2^{7/2 - 1/r} 
\sqrt{\varphi_{r,3,1}(\la_{r,3,1})} }~. \feq
\beq \begin{tabular}{c|c|c|c|c|c} $r$ & 2.2 & 3 & 6 & 10 & 20 \\ \hline
$S^{+}_{r,3,1}$ & 0.8605 & 0.6475 & 0.4888 & 0.4519 & 0.4341 \\ \hline
$S^{-}_{r,3,1}$ & 0.8597 & 0.6458 & 0.4872 & 0.4507 & 0.4333 
\end{tabular}
~~~~. 
\feq
For each $r$ in this table
$(S^{+}_{r,3,1} - S^{-}_{r,3,1})/S^{-}_{r,1,1} < 0.004$, 
with a maximum uncertainty for $r=6$. 
\vskip 0.1cm \noindent
\textbf{C) Case $\mbox{\boldmath $n=2, d=2$}$}~. Eq.s \rref{srnd} 
\rref{srid} give $S^{+}_{r,2,2}$ for all $r$, and in particular
\beq S^{+}_{2, 2, 2} = 1~, \qquad S^{+}_{\infty, 2, 2} = 
1 / \sqrt{4 \pi} \simeq 0.2821~.\feq 
The function $\varphi_{r,2,2}$ computed via Eq.\rref{elem}, its 
minimum point $\la_{r,2,2}$ and the integral $I_{r,2,2}$,  
defined by \rref{irnd}, are given by
$$ \varphi_{r,2,2}(\la) = {\la^4 + \la^2 + 1 \over 6~ \la^{2 - 4/r}}, 
~\la_{r,2,2} = \sqrt{ {- 1/r + \sqrt{1 - 3/r^2} \over 1 + 2/r}};~$$
\beq  I_{r,2,2} = \int_{0}^{+\infty} d t~t \left( t K_1(t) \right)^r~.
\feq
The above integral must be computed numerically. 
The final expression for the lower bounds, and some numerical values for
them and for the upper bounds \rref{srnd} are
\beq S^{-}_{r,2,2} = { I_{r,2,2}^{1/r} \over 2^{3/2 - 1/r} \pi^{1/2 -
1/r}
\sqrt{\varphi_{r,2,2}(\la_{r,2,2})} }~, \feq
\beq \begin{tabular}{c|c|c|c|c|c|c} $r$ & 2.1 & 3 & 6 & 18 & 50 & 100 \\
\hline
$S^{+}_{r,2,2}$ & 0.8494 & 0.4557 & 0.2949 & 0.2644 & 0.2694 & 0.2737 
\\ \hline
$S^{-}_{r,2,2}$ & 0.8465 & 0.4455 & 0.2854 & 0.2582 & 0.2659 & 0.2715  
\end{tabular}~~~~. \feq
It is $(S^{+}_{r,2,2} - S^{-}_{r,2,2})/S^{-}_{r,2,2} < 0.04$ for 
all $r$ in this table, with a maximum uncertainty for $r=6$. 
\vskip 0.1cm \noindent
\textbf{D) Case $\mbox{\boldmath $n=2, d=3$}$}~. Eq.s \rref{srnd} 
\rref{srid} give $S^{+}_{r, 2, 3}$ for all $r$, and in particular
\beq S^{+}_{2, 2, 3} = 1~, \qquad S^{+}_{\infty, 2, 3} = 
1/\sqrt{8 \pi} \simeq 0.1995~.\feq 
The function $\varphi_{r,2,3}$ computed from Eq.\rref{remar}, 
its minimum point $\la_{r,2,3}$ and the integral 
$I_{r,2,3}$ provided by \rref{elem} are
$$ \varphi_{r,2,3}(\la) = {\pi (5 \la^4 + 2 \la^2 + 1) 
\over 32~ \la^{3 - 6/r}},~\la_{r,2,3} = \sqrt{
{1 - 6/r + 4 \sqrt{1 + 3/ r - 9/r^2} \over 5 (1 + 6/r)} } ; $$
\beq I_{r,2,3} = \left( {\pi \over 2} \right)^{r/2}~ {2 \over r^3}. \feq
The final expression for the lower bounds, and some numerical values
for them and for the upper bounds are
\beq S^{-}_{r,2,3} = { \pi^{1/r} E(1/r)^3 \over 2^{5/2 - 3/r} 
\sqrt{\varphi_{r,2,3}(\la_{r,2,3})} }~. \feq
\beq \begin{tabular}{c|c|c|c|c|c|c|c|c} $r$ & 2.1 & 3 & 4 & 7 & 11 & 20
& 
100 & 1000 \\ \hline
$S^{+}_{r,2,3}$ & 0.7830 & 0.3118 & 0.2183 & 0.1657 & 0.1594 
& 0.1647 & 0.1864 & 0.1975 \\ \hline
$S^{-}_{r,2,3}$ & 0.7762 & 0.2912 & 0.1986 & 0.1486 & 0.1437 
& 0.1511 & 0.1795 & 0.1960 
\end{tabular} ~~~~. \feq
For these and other values of $r$ in $(2, \infty)$, we always found
$(S^{+}_{r,2,3} - S^{-}_{r,2,3})/S^{-}_{r,2,3} < 0.12$~,
the maximum uncertainty occurring for $r \simeq 7$.  
\vfill \eject 
\end{document}